\documentclass[10pt]{elsart3-1}

\usepackage{amssymb}
\usepackage[latin1]{inputenc}
\usepackage[T1]{fontenc}
\usepackage[english,francais]{babel}

\newtheorem{theorem}{Theorem}[section]
\newtheorem{lemma}[theorem]{Lemma}
\newtheorem{e-proposition}[theorem]{Proposition}

\newtheorem{e-definition}[theorem]{Definition\rm}

\newtheorem{theoreme}{Th\'eor\`eme}[section]

\newtheorem{proposition}[theoreme]{Proposition}

\newtheorem{remarque}{\it Remarque}

\renewcommand{\theequation}{\arabic{equation}}
\def\og{\leavevmode\raise.3ex\hbox{$\scriptscriptstyle\langle\!\langle$~}}
\def\fg{\leavevmode\raise.3ex\hbox{~$\!\scriptscriptstyle\,\rangle\!\rangle$}}

\journal{the Acad\'emie des sciences}
\begin{document}
\vspace*{-1cm}

\begin{frontmatter}




%
\selectlanguage{francais}
\title{In\'egalit\'es de Poincar\'e  cin\'etiques.}



\author[authorlabel1]{Pascal AZERAD},
\ead{azerad@math.univ-montp2.fr}
\author[authorlabel2]{St\'ephane BRULL}
\ead{stephane.brull@math.univ-toulouse.fr}

\address[authorlabel1]{I3M, UMR 5149, Universit\'e Montpellier 2, F-34095 Montpellier cedex 5}
\address[authorlabel2]{IMT, UMR 5219, Universit\'e Paul Sabatier, F-31062 Toulouse cedex 9 }


\selectlanguage{francais}

\begin{abstract}
\noindent 
Dans cette note, nous \'etablissons des in\'egalit\'es de type Poincar\'e pour une  famille d'\'equations cin\'etiques. Nous appliquons ensuite cette in\'egalit\'e au traitement variationnel d'un mod\`ele cin\'etique lin\'eaire. 
\selectlanguage{francais}

\selectlanguage{english}
\noindent{\bf Abstract}
\noindent
\\ \noindent In this note we prove  Poincar\'e type inequalities for 
a family of kinetic equations. We apply this inequality to the variational solution of a linear kinetic model.
\end{abstract}
\end{frontmatter}

\selectlanguage{english}
\section*{Abridged English version}
Let $T>0$ and   $\Omega$ a regular domain of   $\mathbb{R}^d$, not necessarily  bounded. Denote
\begin{eqnarray*}
a=(1,v) \in \mathbb{R}_t \times \mathbb{R}_x^d,\quad \nabla_{t,x} = (\partial_t,\nabla_x), 
\quad \mathcal{R} = (0,T)\times \Omega. 
\end{eqnarray*}
 Consider the ''Lie-Sobolev'' spaces: 
\begin{eqnarray*}
 H(a,\mathcal{R}) = \{ f \in L^2(\mathcal{R});\;  
a \cdot \nabla_{t,x} f \in L^2(\mathcal{R}) \},
\end{eqnarray*}
 equipped with the norms
$
\| f \|_{H(a,\mathcal{R})} = \| f \|_{L_{t,x}^2 (\mathcal{R})} + 
\| a \cdot \nabla_{t,x} f \|_{L_{t,x}^2 (\mathcal{R})}.
$\\
Define $\partial\mathcal{R}^- = (\{0 \}\times \Omega) \bigcup ((0,T)\times \Gamma^-_v)
$ where $(0,T)\times \Gamma^-_v = \{(t,x) \in \mathbb{R}_t \times \mathbb{R}_x^3
;\; v \cdot n_x <0 \}$, $n_x$ being the outer unit normal to $\partial\Omega$.
We prove the following results.
\begin{proposition}\label{kinpoin} 
(Kinetic Poincar\'e inequality.) Let $T >0$ and  $v \in \mathbb{R}^{d} $.
Let $\mathcal{R}$ a regular enough domain of $(0,T)\times\mathbb{R}^{d} $.  Let $f:= f(t,x) \in H_0(a,\mathcal{R},\partial\mathcal{R}^-)$. Then

\begin{eqnarray} \label{kineg}
 \| f \|_{L_{t,x}^{2}}  \leq\; 2T   \;\|\frac{\partial f}{\partial t} + v \cdot \nabla_{x} f  \|_{L_{t,x}^{2}}.
\end{eqnarray}

\end{proposition}

\begin{proposition} \label{inegpoincarevla}
(Vlasov Poincar\'e inequality.) Let $T >0$ and $\mathcal{R}$ a regular enough domain of  $(0,T)\times \mathbb{R}_x^{3}\times  \mathbb{R}_v^{3}$. Let $f:= f(t,x,v) \in H_0(a,\mathcal{R},\partial\mathcal{R}^-)$. 
Then
\begin{eqnarray} \label{inegvla}
\| f \|_{L_{t,x,v}^{2}}  \leq 2T   \;\|\frac{\partial f}{\partial t} + v \cdot \nabla_{x} f + 
(E + v \times B ) \cdot \nabla_v f  \|_{L_{t,x,v}^{2}}.
\end{eqnarray}

\end{proposition}
We apply inequality (\ref{kineg})  to the variational solution by Lax-Milgram lemma of the following linear kinetic model. 
\begin{eqnarray} \label{eeq}
 (\frac{\partial}{\partial t} + v \cdot \nabla_{x}) u &=& G(t,x,v),\hspace*{0.1 in}  t \in ]0,T[,\;x \in \Omega, \; v \in \mathbb{R}^{d},
\\ \label{eeq1}  u(t=0,x,v) &=& u_{0}(x, v),\hspace*{0.1 in} x\in \Omega ,\; v \in \mathbb{R}^{d},
\\ \label{eeq2}  u(t,\sigma ,v) &=& u_{b}( t, \sigma ,v ), \hspace*{0.1 in} \sigma \in \Gamma^{-}_v, \; v \in \mathbb{R}^{d}.
\end{eqnarray}
With a suitable lifting of the boundary and initial conditions (see Lemma \ref{lem})
 problem (\ref{eeq}, \ref{eeq1}, \ref{eeq2}) is reformulated for each $v
\in \mathbb{R}^{d}:$ 
\begin{eqnarray}
  a \cdot \nabla_{t,x} f &=& G(\cdot,\cdot, v) \hspace*{0.1 in} \mbox{in} \hspace*{0.1 in}  \mathcal{R} = (0,T)\times \Omega,\label{eeqdir}
\\  f(t,\sigma ) &=& 0 \hspace*{0.1 in} \mbox{ on } \hspace*{0.1 in} \partial \mathcal{R}^{-} 
=  \left(\{ 0 \}\times\Omega\right) \cup  \left(]0,T[\times\Gamma^{-}_v\right),\label{edir}
\end{eqnarray}
and we obtain the last result
\begin{proposition} \label{poincare}Let $\Omega$ a  domain  of $\mathbb{R}^{d}.$ 
Problem (\ref{eeqdir}, \ref{edir}) has  a unique strong solution $f:=
f(t,x,v)$ such that  $\forall \; v\in \mathbb{R}^{3}$, we have $ \| f(\cdot,\cdot,v) \|_{L_{t,x}^{2}} \leq 2 T\, \| \, G(\cdot,\cdot,v) \, \|_{L_{t,x}^{2}}.$
\end{proposition}
\selectlanguage{francais}
\section{Introduction}
\label{intro}
Le but de cette note est d'\'etablir des in\'egalit\'es de type Poincar\'e dans un cadre cin\'etique. 
Nous appliquons cette in\'egalite \`a la r\'esolution variationnelle de l' \'equation de transport cin\'etique.
\section{In\'egalit\'e de Poincar\'e cin\'etique.}
Soit $T>0$ et  $\Omega$ un ouvert r\'egulier de   $\mathbb{R}^d$, pas n\'ecessairement born\'e. Notons
\begin{eqnarray*}
a=(1,v) \in \mathbb{R}_t \times \mathbb{R}_x^d,\quad \nabla_{t,x} = (\partial_t,\nabla_x), 
\quad \mathcal{R} = (0,T)\times \Omega. 
\end{eqnarray*}
Consid\'erons les espaces de \og Lie-Sobolev\fg{}, espaces naturels introduits dans  \cite{ces}: 
\begin{eqnarray*}
 H(a,\mathcal{R}) = \{ f \in L^2(\mathcal{R});\;  
a \cdot \nabla_{t,x} f  = \frac{\partial f}{\partial t} + v \cdot \nabla_{x} f \in L^2(\mathcal{R}) \},
\end{eqnarray*}
 munis des normes
\begin{eqnarray*}
\| f \|_{H(a,\mathcal{R})} = \| f \|_{L_{t,x}^2 (\mathcal{R})} + 
\| a \cdot \nabla_{t,x} f \|_{L_{t,x}^2 (\mathcal{R})}.
 \end{eqnarray*}
 Le vecteur $a = (1, v)$ \'etant \emph{constant}, on montre que
$H(a,\mathcal{R})$ est un espace de Hilbert et que les fonctions 
$\mathcal{C}^\infty$ sont denses dans $H(a,\mathcal{R})$.
Comme ${\rm div(a)} = 0$, 
$f\in H(a,\mathcal{R})$ implique que 
\begin{eqnarray*}
f\cdot a \in  H(div,\mathcal{R})=
\{ \varphi \in L^2(\mathcal{R});\; 
 \nabla_{t,x} \cdot \varphi \in L^2(\mathcal{R}) \}.
\end{eqnarray*}
On peut donc définir l'op\'erateur de trace normal
\begin{eqnarray*}
\gamma_n : f \in H(a,\mathcal{R}) \mapsto  f\, (a\cdot n) \in H^{-1/2}(\partial\mathcal{R})
\end{eqnarray*}
 o\`u $n= (n_t,n_x) \in  \mathbb{R}_t \times \mathbb{R}_x^{d}$ est le vecteur normal ext\'erieur et $\partial\mathcal{R} = \left( ]0,T[ \times \partial \Omega \right) \cup 
\left( \{0,T\}\times \Omega \right).$
On introduit alors la frontière entrante (resp.  sortante) spatio-temporelle 
$\partial\mathcal{R}^- =\{(t,\sigma)\in \partial\mathcal{R}; \; 
n_t + v\cdot n_x <0 \} $
(resp $\partial\mathcal{R}^+ =\{(t,\sigma)\in \partial\mathcal{R}; \; 
n_t + v\cdot n_x >0 \}$).
\begin{remarque} La frontière entrante regroupe la condition initiale et la condition limite entrante : 
$\partial\mathcal{R}^- = (\{0 \}\times \Omega) \bigcup ((0,T)\times \Gamma^-_v).$  
\end{remarque}
On peut alors d\'efinir l'espace
\begin{eqnarray*} 
H_0(a,\mathcal{R},\partial\mathcal{R}^-) = \{f \in H(a,\mathcal{R});\; 
\gamma_n f  = 0 \; \mbox{sur}
\;\partial\mathcal{R}^- \}.
\end{eqnarray*}
Le principal r\'esultat est  l'in\'egalit\'e de Poincar\'e cin\'etique suivante.

\begin{proposition} \label{inegpoincare}
(In\'egalit\'e de Poincar\'e cin\'etique.) Soit $T >0$ et $v \in \mathbb{R}^{d} $.
Soit $\mathcal{R}$ un ouvert r\'egulier de $(0,T)\times\mathbb{R}^{d}$.  Soit $f:= f(t,x) \in H_0(a,\mathcal{R},\partial\mathcal{R}^-)$. 
On a l'in\'egalit\'e

\begin{eqnarray} \label{ineg}
 \| f \|_{L_{t,x}^{2}}  \leq\; 2T   \;\|\frac{\partial f}{\partial t} + v \cdot \nabla_{x} f  \|_{L_{t,x}^{2}}  .
\end{eqnarray}

\end{proposition}
 (Propri\'et\'e \ref{inegpoincare}). 
On raisonne par densité en prenant $f$ lisse.
Consid\'erons $w(t,x): = t-T.$  Pour tout $(x,t) \in \mathcal{R}, \;
w(t,x)\leq 0$ et
\begin{eqnarray*}
\frac{\partial w}{\partial t} + v \cdot \nabla_{x} w = 1.
\end{eqnarray*}
D'apr\`es la formule de Stokes, on obtient
\begin{eqnarray*}
\int \int_\mathcal{R} \left(\frac{\partial }{\partial t} + v \cdot \nabla_{x}\right) (w\cdot f^2)
= \int_{ \partial \mathcal{R}}  w f^2\, (n_t + v\cdot n_x).
\end{eqnarray*}
Or $f(n_t + v\cdot n_x)=\gamma_n f= 0$  sur $\partial \mathcal{R}^-$. Donc comme $w\leq 0$ et $f^2 \geq 0$, on a
\begin{eqnarray*}
\int \int_\mathcal{R} \left(\frac{\partial }{\partial t} + v \cdot \nabla_{x}\right) (w\cdot f^2) 
=  \int_{ \partial\mathcal{R}^+}  w f^2\, (n_t + v\cdot n_x)\leq 0.
\end{eqnarray*}
D'autre part,
\begin{eqnarray*}
\left(\frac{\partial }{\partial t} + v \cdot \nabla_{x}\right) (w\cdot f^2)
&=&  f^2 + 2 w\, f\, (\frac{\partial f}{\partial t} + v \cdot \nabla_{x} f).
\end{eqnarray*}
On obtient alors
\begin{eqnarray*}
\int \int_Q   f^2 \leq 
2 \|w\|_{L^\infty} \int \int_Q \left|f\right| 
\left|\frac{\partial f}{\partial t} + v \cdot \nabla_{x} f\right|.
\end{eqnarray*}
Avec  l'in\'egalit\'e de Cauchy-Schwarz 
\begin{eqnarray*}
\| f \|_{L_{t,x}^{2}}  \leq 2 \|w\|_{L^\infty}\;\|\frac{\partial f}{\partial t} + v \cdot \nabla_{x} f  \|_{L_{t,x}^{2}}.
\end{eqnarray*}
Comme $\|w\|_{L^\infty} \leq T $, on obtient le r\'esultat.
\begin{remarque}
La constante $T$ est ind\'ependante de $v$ et le domaine $\Omega$ n'a pas besoin d'\^etre born\'e. 
En effet le domaine espace-temps $\mathcal{R}$ est automatiquement born\'e dans la direction du \emph{temps}.
\end{remarque}
\vspace{-0.7cm}
\section{Une in\'egalit\'e de  Vlasov-Poincar\'e.}
A partir de maintenant $d=3$. Soit $\mathcal{R} = (0,T)\times \Omega \times \mathbb{R}^3$. Consid\`erons l'\'equation de Vlasov 
\begin{eqnarray*}
\frac{\partial f}{\partial t } + v \cdot \nabla_x f + ( E + v \times B ) \cdot \nabla_v f =0
\end{eqnarray*}
o\`u $E: = E(t,x)$ et $B:= B(t,x)$ sont respectivement les champs \'electriques et magn\'etiques, suppos\'es r\'eguliers, par exemple $\mathcal{C}^1$. 
 Notons
$$a=\left(1, v, E(t,x) + v \times B(t,x)\right) \in
 \mathbb{R}_t \times \mathbb{R}_x^3 \times \mathbb{R}_v^3.$$
On remarque que  
$\nabla_{t,x,v}\cdot a = \nabla_v \cdot (v \times B(t,x))=
\sum_j \frac{\partial\mbox{}}{\partial v_j}\;(v\times B(t,x))_j = 0,
$
car $(v \times B(t,x))_j$ d\'epend seulement de $v_i$ pour 
$i\neq j$. 
On a alors
\begin{equation}
\nabla_{t,x,v}\cdot (f \cdot a) = a \cdot \nabla_{t,x,v} f
\label{solenoid}  
\end{equation}
\begin{remarque}
La propri\'et\'e (\ref{solenoid}) permettrait certainement d'abaisser la r\'egularit\'e de $E$ et $B$.
\end{remarque}
On consid\`ere l'espace de Sobolev anisotrope: 
\begin{eqnarray*}
 H(a,\mathcal{R}) = \{ f \in L^2(\mathcal{R});\; f \cdot a \in L^2(\mathcal{R}) \;
\mbox{et} \;  
a \cdot \nabla_{t,x,v} f \in L^2(\mathcal{R}) \},
\end{eqnarray*}
muni de la norme
\begin{eqnarray*}
\| f \|_{H(a,\mathcal{R})} = \| f\cdot a \|_{L_{t,x,v}^2 (\mathcal{R})} + 
\| a \cdot \nabla_{t,x,v} f \|_{L_{t,x,v}^2 (\mathcal{R})}.
 \end{eqnarray*}
 \begin{remarque} 
Comme $v$ peut \^etre arbitraire, $a$ est non born\'e. 
 \end{remarque}
D'apr\`es (\ref{solenoid}) on a
\begin{eqnarray*}
 H(a,\mathcal{R}) = \{ f \in L^2(\mathcal{R});\; f \cdot a \in H(div,\mathcal{R}) \}.
\end{eqnarray*} 
On peut alors d\'efinir un op\'erateur de trace normal
\begin{eqnarray*}
\gamma_n : f \in H(a,\mathcal{R}) \mapsto  f\, (a\cdot n) \in H^{-1/2}(\partial\mathcal{R})
\end{eqnarray*}
 o\`u $n= (n_t,n_x,0) \in  \mathbb{R}_t \times \mathbb{R}_x^{3} \times \mathbb{R}_v^{3} $ est la normale extérieure au bord
$$\partial\mathcal{R} = \left( ]0,T[ \times \partial \Omega 
\times \mathbb{R}_v^{3}  \right) \cup 
\left( \{0,T\}\times \Omega \times \mathbb{R}_v^{3} \right).$$
Notons qu'il n'y a pas de bord dans la direction $v$.
On a $a\cdot n = n_t + v\cdot n_x$ et le bord entrant est 
$$\partial\mathcal{R}^- = \left( ]0,T[ \times \Gamma^-_v 
\times \mathbb{R}_v^{3}  \right) \cup 
\left( \{0\}\times \Omega \times \mathbb{R}_v^{3} \right).$$
On d\'efinit alors l'espace
\begin{eqnarray*} 
H_0(a,\mathcal{R},\partial\mathcal{R}^-) = \{f \in H(a,\mathcal{R});\;  f (a\cdot n) = 0 \; \mbox{sur}
\;\partial\mathcal{R}^- \}.
\end{eqnarray*} La proposition s'\'enonce.
\begin{proposition} \label{inegpoincarevla}
(Inégalité de Vlasov Poincar\'e.) Soit $T >0$ et $\mathcal{R}$ un ouvert r\'egulier de $(0,T)\times \mathbb{R}_x^{3}\times  \mathbb{R}_v^{3}$. Let $f:= f(t,x,v) \in H_0(a,\mathcal{R},\partial\mathcal{R}^-)$. 
On a l'in\'egalit\'e

\begin{eqnarray} \label{inegvla}
\| f \|_{L_{t,x,v}^{2}}  \leq 2T   \;\|\frac{\partial f}{\partial t} + v \cdot \nabla_{x} f + 
(E + v \times B ) \cdot \nabla_v f  \|_{L_{t,x,v}^{2}}  .
\end{eqnarray}

\end{proposition}
Preuve (analogue \`a la proposition \ref{ineg}). On consid\`ere $w(t,x): = t-T.$ Pour tout $(x,t) \in \mathcal{R} $, 
$w(t,x)\leq 0$ et 
\begin{eqnarray*}
\frac{\partial w}{\partial t} + v \cdot \nabla_{x} w + (E + v \times B ) \cdot \nabla_v w  = 1.
\end{eqnarray*}
D'apr\`es la formule de Stokes, on obtient
\begin{eqnarray*}
\int\int \int_\mathcal{R} \left(\frac{\partial }{\partial t} + v \cdot \nabla_{x} +
( E + v \times B) \cdot \nabla_v \right) (w\cdot f^2)(t,x,v) \, dt dx dv
\end{eqnarray*}
\begin{eqnarray*}
= \int\int \int_\mathcal{R} \nabla_{t,x,v}  (w f^2 \cdot a)
=\int \int_{ \partial \mathcal{R}}  w f^2\, (n_t + v\cdot n_x)
= \int \int_{ \partial\mathcal{R}^+}  w f^2\, (n_t + v\cdot n_x)\leq 0.
\end{eqnarray*}
Or 
$
\left(\frac{\partial }{\partial t} + v \cdot \nabla_{x} +
( E + v \times B) \cdot \nabla_v \right) (w\cdot f^2)  =  f^2 + 2\, w \,f \left(\frac{\partial f}{\partial t} 
+ v \cdot \nabla_{x} f+
( E + v \times B) \cdot \nabla_v f \right),
$
donc
\begin{eqnarray*}
 \int_{\mathcal{R}}   f^2 \leq 
2 \|w\|_{L^\infty}  \int_{\mathcal{R}}| f |
\left|\frac{\partial f}{\partial t} + v \cdot \nabla_{x} f +
( E + v \times B) \cdot \nabla_v f\right|
\end{eqnarray*}
et on conclut comme dans la preuve de (\ref{inegpoincare}) $\Box$
\section{Application \`a la forme variationnelle du transport cin\'etique.}
\noindent Soit $v \in \mathbb{R}^{3} \mapsto  G(\cdot,\cdot,v) \in L^2((0,T)\times \Omega)$. Consid\'erons le probl\`eme suivant, pour chaque $v\in \mathbb{R}^{3}$.
\begin{eqnarray} \label{eq}
 (\frac{\partial}{\partial t} + v \cdot \nabla_{x}) u &=& G(t,x,v),\hspace*{0.1 in}  t \in ]0,T[,\;x \in \Omega,
\\ \label{eq1}  u(t=0,x,v) &=& u_{0}(x, v),\hspace*{0.1 in} x\in \Omega ,
\\ \label{eq2}  u(t,\sigma ,v) &=& u_{b}( t, \sigma ,v ), \hspace*{0.1 in} \sigma \in \Gamma^{-}_v.
\end{eqnarray}
 Le probl\`eme (\ref{eq}, \ref{eq1}, \ref{eq2}) est reformul\'e  pour chaque $v\in \mathbb{R}^{3}$:
\begin{eqnarray*}
  a \cdot \nabla_{t,x} u &=& G(\cdot,\cdot, v) \hspace*{0.1 in} \mbox{dans} \hspace*{0.1 in}  \mathcal{R} = (0,T)\times \Omega,
\\  u(t,\sigma ) &=& g(t, \sigma,v ) \hspace*{0.1 in} \mbox{ sur } \hspace*{0.1 in} \partial \mathcal{R}^{-} 
=  \left(\{ 0 \}\times\Omega\right) \cup  \left(]0,T[\times\Gamma^{-}_v\right).
\end{eqnarray*}
Pour se ramener à des conditions de bord de type Dirichlet 
homog\`ene on d\'emontre  ais\'ement le lemme suivant par la m\'ethode des caract\'eristiques, qui sont ici rectilignes.
\begin{lemma} \label{lem}
Soit $\Omega$ un domaine de $\mathbb{R}^{3}.$ 
Le vecteur $v \in \mathbb{R}^{3}$ \'etant donn\'e, il existe $g(t,x,v)$ solution du probl\`eme
\begin{eqnarray*}
  a \cdot\nabla_{t,x} g &=&0,\hspace*{0.1 in}  t \in ]0,T[,\; x \in \Omega
\\ g(0,x,v) &=& u_{0}(x,v), \hspace*{0.1 in} x \in \Omega,
\\ g(t,\sigma ,v) &=& u_{b}(t,\sigma ,v ), \hspace*{0.1 in}  t \in ]0,T[,\; \sigma \in \Gamma^{-}_v.
\end{eqnarray*}
De plus, si $u_{0}(\cdot,v) \in L^{\infty}(\Omega )$ et 
si $u_{b}(\cdot,\cdot, v) \in L^{\infty}((0,T) \times  \Gamma^- _v)$ alors 
$g(\cdot,\cdot,v) \in L^{\infty} ([0,T] \times \Omega )$ et
\begin{eqnarray*}
 \|g(\cdot,\cdot,v)\|_{L^{\infty}([0,T] \times \Omega)}\leq \| u_{0}(\cdot,v) \|_{L^{\infty}(\Omega)} + \| u_{b}(\cdot,\cdot,v) \|_{L^{\infty}([0,T] \times \partial \Omega)}.
 \end{eqnarray*}
\end{lemma}
On effectue alors le changement d'inconnue
 $f = u-g$ o\`u $g$ est le rel\`evement d\'efini par le lemme \ref{lem}. Ainsi $f$ devient solution du probl\`eme de Dirichlet homogène suivant
\begin{eqnarray}
  \label{eqdir} a \cdot \nabla_{t,x} f(t,x,v) &=& G(t,x,v) \hspace*{0.1 in} \mbox{dans} \hspace*{0.1 in}  \mathcal{R},
\\  \label{dir} f(t,\sigma ,v ) &=& 0 \hspace*{0.1 in} \mbox{sur} \hspace*{0.1 in}  \partial\mathcal{R}^-.
\end{eqnarray}
\noindent  On obtient alors la proposition suivante.
\begin{proposition} \label{poincare}Soit $\Omega$ un domaine  de $\mathbb{R}^{3}.$ 
Le probl\`eme (\ref{eqdir}, \ref{dir}) poss\`ede une unique solution forte $f:=
f(t,x,v)$ v\'erifiant pour chaque $v\in \mathbb{R}^{3}$
\begin{eqnarray} \label{ineg}
 \| f(\cdot,\cdot,v) \|_{L_{t,x}^{2}} \leq C\, \| \, G(\cdot,\cdot,v) \, \|_{L_{t,x}^{2}}
\end{eqnarray}
o\`u $C$ est une constante positive ind\'ependante de la variable $v$ et born\'ee par $2T$.
\end{proposition}
\vspace*{1 mm}
\begin{remarque}
La preuve est inspir\'ee de la m\'ethode STILS d\'evelopp\'ee dans
 (\cite{AP}, \cite{APe}, \cite{BP}),  adapt\'ee \`a un cadre cin\'etique. 
 Mais la constante $C$ est ind\'ependante de la variable $v$ contrairement \`a (\cite{AP}).    
\end{remarque}
\vspace*{1 mm}
 (Proposition \ref{poincare}).
Pour chaque $v$ fix\'e, consid\'erons la forme bilin\'eaire suivante 
\begin{eqnarray*}
 \mathcal{B}(f,g) = \int\int_{\mathcal{R}}(a\cdot \nabla_{t,x}\, f)\, (a\cdot \nabla_{t,x}\, g)\,dtdx
 \end{eqnarray*}
et la forme lin\'eaire
\begin{eqnarray*}
L(g) = \int\int_{\mathcal{R}} G \,( a\cdot \nabla_{t,x}\,g).
\end{eqnarray*}
La formulation STILS (\cite{APe,AP}) du probl\`eme s'\'ecrit:
 $G \in L^2(\mathcal{R})$ \'etant donn\'e, trouver
$f \in H_0(a,\mathcal{R}, \partial \mathcal{R}^- ) $ telle que 
\begin{equation}
\mathcal{B}(f,g) = L(g),\quad \forall g \in H_0(a,\mathcal{R}, \partial \mathcal{R}^-).
\label{stils}
\end{equation}
D'apr\`es l'in\'egalit\'e de Poincar\'e \ref{inegpoincare}, la forme bilin\'eaire $\mathcal{B}$ est coercive.
En utilisant le th\'eor\`eme de Lax-Milgam, on obtient  que le probl\`eme (\ref{stils}) est bien pos\'e.
Il reste \`a prouver que $f$ satisfait \`a 
(\ref{eqdir}, \ref{dir}) fortement. Il suffit alors de prouver que 
\begin{eqnarray*}
\{\frac{\partial \varphi}{\partial t} + v \cdot \nabla_{x} \varphi \in L^2_{t,x}; \; \varphi \in H_0(a,\mathcal{R}, \partial \mathcal{R}^-)\}
\end{eqnarray*}
 est dense dans $L^2_{t,x}$. 
Pour cela on reproduit la preuve donn\'ee dans \cite{Az} 
(Th\'eor\`eme 16, pp. 83-84). $\Box$

\end{document}